\documentclass{siamltex}

\usepackage{epsfig}
\usepackage{amsfonts}

\title{Transforming Triangulations on Non-Planar Surfaces}

\author{C. Cort\'es\footnotemark[2] \and C.~I. Grima\footnotemark[2] \and F. Hurtado\footnotemark[3] \and A. M\'arquez\footnotemark[2] \and F. Santos\footnotemark[4] \and J. Valenzuela\footnotemark[2]}

\begin{document}

\maketitle

\renewcommand{\thefootnote}{\fnsymbol{footnote}}

\footnotetext[2]{Dept. Matem\'atica Aplicada I, Univ. de Sevilla, Spain (\{ccortes, grima, almar, jesusv\}@us.es). Partially supported by Projects
MTM2008-05866-C03-01 and P06-FQM-01649}
\footnotetext[3]{Dept. Matem\'{a}tica Aplicada I, Univ. Polit\'{e}cnica de Catalunya, Spain (ferran.hurtado@upc.edu). Partially
supported by Projects MICINN MTM2009-07242 and Gen. Cat.  2009SGR1040}
\footnotetext[4]{Dept. Matem\'aticas, Estad\'{\i}stica y Computaci\'on, Univ. de
Cantabria, Spain. (santosf@unican.es). Partially supported by MTM2008-04699-C03-02}

\renewcommand{\thefootnote}{\arabic{footnote}}

\begin{abstract}
We consider whether any two triangulations of a polygon or a point set on a non-planar surface
with a given metric can be transformed into each other by a sequence of edge flips. The answer is
negative in general with some remarkable exceptions, such as polygons on the cylinder, and on the
flat torus, and certain configurations of points on the cylinder.
\end{abstract}

\begin{keywords}
Graph of triangulations, triangulations on surfaces, triangulations of
polygons, edge flip.
\end{keywords}

\begin{AMS}
68U05, 52C99, 65D18, 68R10
\end{AMS}

\pagestyle{myheadings}
\thispagestyle{plain}
\markboth{CORT\'{E}S, GRIMA, HURTADO, M\'{A}RQUEZ, SANTOS AND VALENZUELA}{TRANSFORMING TRIANGULATIONS ON NON-PLANAR SURFACES}

\section{Introduction}

Most of the problems considered so far in Computational Geometry are restricted to the plane, or
to the Euclidean 3-space. However, in many applications it is necessary to deal with input data
that lies on a surface rather than in the plane. Recently, some works have been focused in solving
some of  the problems arising in those cases (cf. \cite{albertoyclara,mazon,okabe}). This paper is
part of this stream, studying the graph of triangulations of a polygon on a surface.

Partitioning geometric domains into simpler pieces, such as
triangles, is a common strategy  to several fields, the finite
element method being a most relevant example. In particular, the
triangulation of polygons is an intermediate step in many
algorithms in the area of Computational Geometry.

In many cases, it is not only needed to obtain a triangulation of
a given region, but a ``good" one. Some examples of this
assertion can be found when it is desired to improve the quality
of a graphic representation or to find a ``nice" mesh on a given
surface in order to apply finite elements methods. When the {\em
quality\/} of the triangulation with respect to some criterion is
considered, and no direct method for obtaining the optimal
triangulation is known, it is natural to perform operations that
allow local improvements. The best-known method is the edge flip:
when two triangles form a convex quadrilateral, their common edge
is replaced by the other diagonal of the
quadrilateral~\cite{bern,fortune2}. This local transformation,
introduced by Lawson in \cite{lawson2}, can be combined if
necessary with methods such as simulated annealing to escape
local optima ~\cite{gelatt,vanlaar} and has also been used for
the purposes of enumeration~\cite{avis-fukuda}. It also admits
several variations~\cite{pocch,santos:bistellar}. Regarding the
local operation we have just described, a basic issue is whether
any two triangulations of a domain $D$ can be transformed into
each other by means of a sequence of flips. If we define the {\em
triangulation graph} of $D$, as that graph $TG(D)$ having as nodes
the triangulations of $D$, with adjacencies corresponding to edge
flips, the above question becomes obviously whether $TG(D)$ is a
connected graph or not.

It is known that the graph of triangulations of a planar simple polygon or a point set with $n$
vertices is connected and its diameter is $O(n^2)$, which is tight~\cite{hurt2}. It is worth
mentioning that even the case of a convex $n$-gon $P$ has been thoroughly studied because $TG(P)$
is isomorphic to the rotation graph of binary trees with $n-2$ internal nodes~\cite{hurt,sleator}.
On the sphere, the situation is essentially the same as in the plane.

In this work we study the connectivity of the triangulation graph
for simple polygons and point sets lying on surfaces. Regarding
polygons we prove that for the cylinder and the torus with
their flat metrics the graph is always connected (if non-empty). For general surfaces
and metrics the situation is usually the opposite. Even worse, for point sets only certain configurations on the
cylinder have  a connected graph of triangulations.

At this point it is convenient to clarify that with the general
purpose of extending Computational Geometry to surfaces, it is
needed to ``translate'' some of the elements that usually appear
in the plane to the surfaces; in our case, we need to know how to
join a pair of points (in other words, how to translate the
concept of segment); it is known that, in general, there are
infinitely many geodesics joining two points, but usually only
one with the minimal length (see~\cite{docarmo2}). Thus,
following the cited works \cite{albertoyclara,mazon,okabe} and
others, this unique minimal geodesic joining a pair of points on
a surface will be called {\em the segment\/} defined by that pair
of points. In what follows, only segments between pairs of points
will be considered.

Equally some words must be said about the surfaces, or, more concretely about the metric, that we are considering here. In general, we will study the case of the locally Euclidean surfaces (those surfaces isometric to the plane in sufficiently small regions). These surfaces have two advantages, on one hand, they are general enough in order to model many practical cases or approximate some other metrics, and, on the other, they have an easy representation as we will se in the next section. Nevertheless, in \S\ref{ccsurfaces} the results are presented in a more general context because we do not need the flat representation of the locally Euclidean surfaces (although an alternative proof of the main result of this section is presented later in the context of locally Euclidean surfaces).

The paper is organized as follows. In \S\ref{prelim} we give definitions and preliminary
results, and we establish the notation that will be used along this paper.
\S\ref{ccsurfaces} shows one of the main results of this paper, that in every compact
connected surface it is always possible to find a metric that admits polygons and point sets with
non-connected graph of triangulations. \S\ref{locally} focus on the connectivity of the
graph of triangulations for both polygons and point sets on the locally Euclidean surfaces. We
conclude in \S\ref{open-p} with some comments and open problems.

\section{Preliminaries}
\label{prelim}

As it is known, many practical problems cannot be modeled by planar situations and other surfaces
are required. When we meet phenomena in which the same configuration of generating points appears
in cycles, we may analyze them with the aid of a point configuration on the cylinder or the torus.
These are two well known surfaces since, together with the twisted cylinder (or infinite M\"obius
strip) and the Klein bottle, they easily admit quotient metrics that make them locally Euclidean.
With these metrics, the graph of triangulations of a polygon both on the cylinder and on the torus
is connected, although this fact does not hold on the other two non-orientable surfaces.

We start this section summarizing the basic properties of the
locally Euclidean surfaces, via their planar representation. A
more complete study of them can be found in~\cite{nikulin}.

\subsection{Locally Euclidean surfaces}\label{locallyEuclidean}
A 2-dimensional locally Euclidean surface is a surface which is isometric with the plane in
sufficiently small regions.

A {\em motion\/} in the plane is a map that preserves distances between points. The group of
motions in the plane is denoted by Mo$(\mathbb{R}^2)$ and consists of translations, rotations, reflections
and glide reflections.

A group $\Gamma \subseteq$ Mo$(\mathbb{R}^2)$ is said to be {\em uniformly discontinuous\/} if there exists
a positive number $d$ such that if $\gamma$ is a motion in $\Gamma$ and $P$ any point in the plane
being $\gamma(P) \neq P$, then the distance between $P$ and $\gamma(P)$ is greater or equal than
$d$.

There are five different types of uniformly discontinuous groups of motions of the plane, up to
isomorphisms: Types I, II.a, II.b, III.a and III.b  \cite{nikulin}, and they can be generated as
follows:
\begin{itemize}
\item Type I, generated by the identity motion.
\item Type II.a, generated by a translation.
\item Type II.b, generated by a glide reflection.
\item Type III.a, generated by two non-collinear translation vectors.
\item Type III.b, generated by a translation and a glide reflection, the direction
of the translation vector being orthogonal to the axis of the glide reflection.
\end{itemize}

Given a group $\Gamma \subseteq$ Mo$(\mathbb{R}^2)$ and a point $P$ in the plane, the {\em orbit\/} of $P$
via $\Gamma$, denoted $\Gamma (P)$, is the set of the successive images of $P$ under the action of
the elements of $\Gamma$, that is $\Gamma(P)=\{\gamma(P) : \gamma \in \Gamma \}$. For any
uniformly discontinuous group of motions $\Gamma \subseteq$ Mo$(\mathbb{R}^2)$ the following notion of
equivalence on points in the plane can be defined: points $A$ and $B$ are equivalent if they
belong to the same orbit namely, there exists a motion $\gamma \in \Gamma$ such that
$\gamma(A)=B$. The orbits are then the equivalence classes under this relation. The set of all
orbits of $\mathbb{R}^2$ under the action of $\Gamma$ is written as $\mathbb{R}^2 / \Gamma$ and is called the {\em
quotient space}. The distance between two points (orbits) ${\bf A}=\Gamma(A)$ and ${\bf
B}=\Gamma(B)$ in $\mathbb{R}^2 / \Gamma$ is defined to be the shortest of the distances $|AB|$, where $A$
and $B$ are points of the plane with $A$ belonging to ${\bf A}$ and $B$ to ${\bf B}$.

Every locally Euclidean surface $\Sigma$ corresponds to a uniformly discontinuous group $\Gamma$
of motions of the plane, so that $\Sigma$ can be obtained from $\Gamma$ as the quotient space
$\mathbb{R}^2 / \Gamma$. Hence there are exactly five types of locally Euclidean surfaces \cite{nikulin}:
the plane (Type I), the cylinder (Type II.a), the twisted cylinder (Type II.b), the (flat) torus (Type
III.a), and the Klein bottle (Type III.b).  Although the term {\em flat torus} applies to surfaces generated by any group of motions of Type III.b, we will follow the convention that consider the translations to be orthogonal. If the translations are not orthogonal then we call the surface so obtained a {\em skew torus}. As we will see in \S\ref{torus} this distinction is not trivial and has important consequences on the connectivity of the graph of triangulations.

According to the above definitions and results, a point $a$ of the surface defined by a  uniformly
discontinuous group $\Gamma$ is specified by an orbit ${\bf A}$ of $\Gamma$. However, in order to
specify $a$, there is no need to know all points of ${\bf A}$; we need only know one point $A$ of
${\bf A}$, and then all the others are obtained from $A$ by applying motions in the given group
$\Gamma$. Therefore, in order to determine the set of all points of the surface we need only
specify some region of the plane, for example a polygon, satisfying the following properties:

\begin{remunerate}
\item The region contains one point from every set of equivalent
points of the plane.
\item No interior point of the region is equivalent to any other
point of the region; that is, equivalent points of the region can
only lie on the boundary.
\end{remunerate}

A region in the plane satisfying (1) and (2) is called a {\em fundamental domain} and the set of
points of the surface is obtained from this region by identifying or gluing together equivalent
points of its boundary. In general we will use the fundamental domains that are more common in the literature, that is, an infinite band for both the cylinder and the twisted cylinder; and a rectangle on the torus and the Klein bottle. In the skew torus it is also usual to consider as fundamental domain a parallelogram which sides are parallel to the direction of the translations.

\begin{figure}[htb]
\begin{center}
\includegraphics[height=4cm]{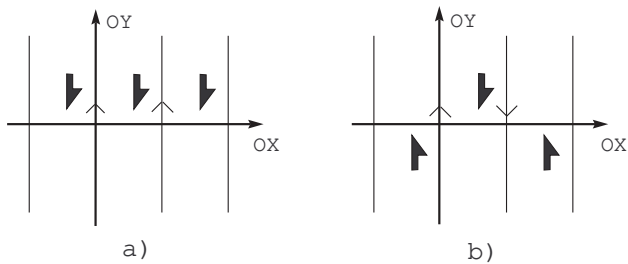}
\caption{The orbit of a polygon in; a) the cylinder, and b) the
twisted cylinder.} \label{cyls}
\end{center}
\end{figure}

\begin{figure}[htb]
\begin{center}
\includegraphics[height=5cm]{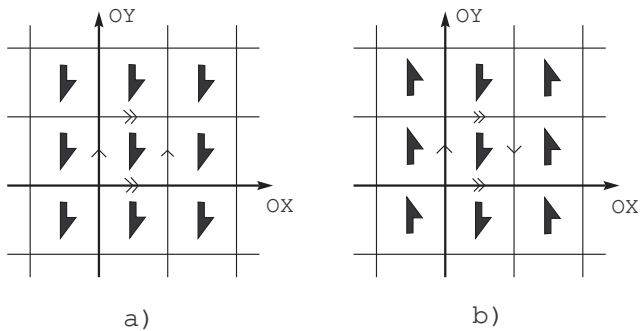}
\caption{The orbit of a polygon in; a) the torus, and b) the Klein
bottle.} \label{tor-kl}
\end{center}
\end{figure}

In order to fix the points in the examples given in \S\ref{locally} we will consider an orthogonal reference system in these surfaces which will be centered, for simplicity, in the leftmost side of the band or in the lowest leftmost corner of the
rectangle (or parallelogram) considered as fundamental domain. In the non-orientable case the $OX$ axis will be taken to coincide with one glide
reflection axis of $\Gamma$. The tesselations of the plane generated by the previous fundamental domains of
each surface together with the orbit of a polygon are depicted in Figures~\ref{cyls} and
\ref{tor-kl}.

\subsection{Triangulations of Euclidean polygons. Flips}\label{trianpolygon}

A {\em Euclidean polygon\/} in a locally Euclidean surface is a region homeomorphic to a closed
disc and whose boundary consists of finitely many geodesic arcs. A Euclidean polygon may be represented
as a simple planar polygon although, depending of the election of the fundamental domain, it might not be completely contained only in one of them.

From now on, Euclidean polygons will be assumed to be already drawn
in the plane.

\begin{figure}[htb]
\begin{center}
\includegraphics[height=4.5cm]{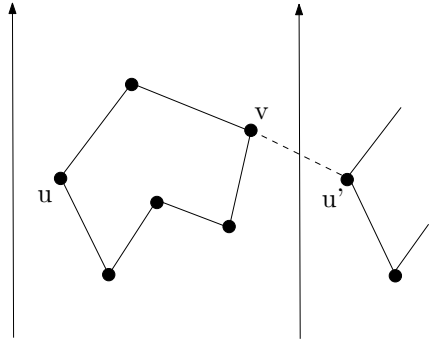}
\caption{Since the nearest copy of $u$ from $v$ is $u'$, $u$ and $v$ cannot be matched inside the
polygon and diagonal $\overline{uv}$ is not admissible.} \label{noninner}
\end{center}
\end{figure}

The segment (that is, the minimum geodesic) between two non-consecutive vertices of a Euclidean
polygon is called a {\em diagonal} of the polygon. The diagonal $\overline{uv}$ is said to be {\em
admissible} if it is contained inside the polygon (Figure~\ref{noninner}).

A {\em (metrical) triangulation\/} of a Euclidean polygon is a partition of the polygon into
triangular regions (that is, regions homeomorphic to a disc bounded by three segments) by means of {\em admissible} diagonals with no intersections except for their ends. Note that we force every face of a triangulation to be triangular instead of considering a
maximal set of segments since, despite being equivalent definitions in the plane, this is no
longer true in other surfaces, as will be apparent in \S\ref{torus}. In the same way, we
define triangulations of point sets as a maximal set of non-crossing segments such that
each bounded region is triangular.
On the contrary what happens in Euclidean polygons, given a point set the shape of the region triangulated depends on the position of the points on the surface, and it can be a Euclidean polygon, or a strip bounded by two geodesics, or the whole surface (see \cite{euclideanposition,albertoyclara}).

\begin{figure}[h]
\begin{center}
\includegraphics[height=6cm]{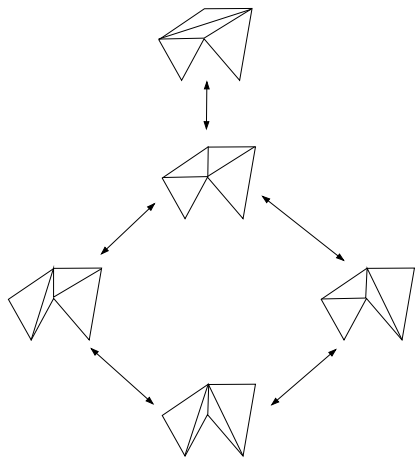}
\caption{The graph of triangulations of a polygon in the plane.}
\label{tri27}
\end{center}
\end{figure}

Let $\{v_i,v_j,v_k\}$ and $\{v_i,v_j,v_l\}$ be two triangles in a triangulation sharing the
diagonal $\overline{v_iv_j}$. By {\em flipping\/} $\overline{v_iv_j}$ we mean the operation of
removing $\overline{v_iv_j}$ and replacing it by the other diagonal $\overline{v_kv_l}$, if it is
admissible in the quadrangle $\{v_i,v_k,v_j,v_l\}$. The {\em graph of triangulations\/} of a
polygon or a point set $P$ is the graph $TG(P)$ having as nodes the triangulations of $P$, with
adjacencies corresponding to diagonal flips (Figure~\ref{tri27}).

\section{Graph of triangulations of a polygon on non-planar surfaces}
\label{ccsurfaces}

One expect that metrical triangulations depend strongly on the metric considered since small
changes in the metric might turn admissible diagonals into non-admissible ones and flip
performance would be affected. In this section, we define a metric on the sphere that produces
polygons and point sets with non-connected graph of triangulations. The same idea will be used to
extend this result to a general closed connected surface.

\begin{figure}[htb]
\begin{center}
\includegraphics[height=3cm]{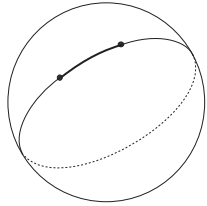}
\caption{The distance between two points in the sphere is given by
the shortest arc of the great circle joining the points.}
\label{sphere}
\end{center}
\end{figure}

On the sphere, with its natural metric, geodesics correspond to great circles and the distance
between two points is the length of the shortest arc of the great circle joining them
(Figure~\ref{sphere}), which is
unique with the exception of {\em antipodal} (or diametrically opposite) points. A (Euclidean) polygon on the sphere, as in a locally Euclidean surface, is
a region homeomorphic to a closed disc and whose boundary consists of finitely many geodesic arcs.
Triangulations, flips and graphs of triangulations of polygons on the sphere are also defined in
the same way as they were in the previous section.

By using similar arguments to those in~\cite{lawson2}, it can be established that the graph of
triangulations of any polygon on the sphere is connected with this metric. But it is possible to
disturb slightly the metric so that this assertion will no longer be true.

\begin{lemma} \label{lsnc}
There exists a surface $M$ homeomorphic to the sphere (in other words, $M$ is a sphere with a metric other than the Euclidean distance) such that in $M$ there exists a Euclidean polygon with a non-connected graph of triangulations and a point set also with a non-connected graph of triangulations.
\end{lemma}

\begin{proof}
Consider a great circle $C$ that divides the sphere into two open hemispheres $H_1$ and $H_2$. Let $p_1,p_2,...,p_6$ be a sequence of vertices uniformly distribute on $C$ such that the great circles joining $(p_1,p_4)$, $(p_2,p_5)$ and $(p_3,p_6)$ intersect only in two antipodal points $n$ and $s$ in $H_1$ and $H_2$ respectively. We move a little bit the vertices $p_1,p_2,...,p_6$ towards $n$ until the arc joining them inside $H_1$ is slightly shorter than the one that crosses through $H_2$.

Let $L=<p_1,p_2,...,p_6>$  be a closed polygonal chain strictly contained in $H_1$, and $P$ the polygon bounded by $L$ which interior is the region with smaller area of the two in which the surface is divided. Now, $M$ is obtained from the
sphere by lifting up a small region around $n$ until the distances (considering the metric inherit from $\mathbb{R}^3$) between
$(p_1,p_4)$, $(p_2,p_5)$ and $(p_3,p_6)$, are enlarged enough to ensure that the diagonals joining
them are non-admissible in $P$ (so those admissible diagonals are exterior to $P$), but without changing the length of the other diagonals of $P$ (Figure~\ref{pol-sph} a)).

\begin{figure}[htb]
\begin{center}
\includegraphics[height=4.5cm]{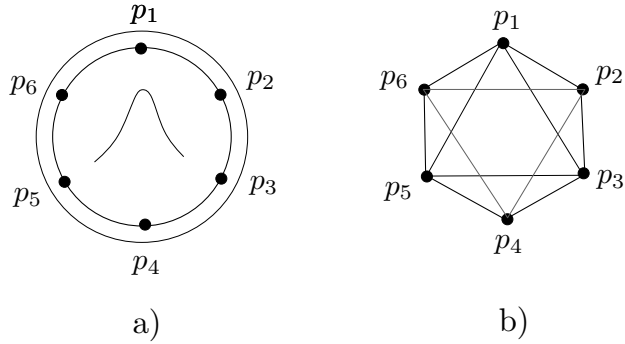}
\caption{An hexagon with two disjoint triangulations in a ``mountainous" sphere.} \label{pol-sph}
\end{center}
\end{figure}

After the lifting of the region around $n$, the length of any geodesic inside $H_1$ on $M$ either is increased or remains the same that its length before the lifting. Moreover, the segments (shortest geodesic arcs) joining $(p_1,p_4)$, $(p_2,p_5)$ and $(p_3,p_6)$ are the arcs of the great circles that join those points in $H_2$.
Therefore, $P$ admits only two different triangulations, shown in Figure~\ref{pol-sph} b), which
cannot be transformed into each other by a sequence of flips, and hence, the graph of
triangulations of $P$ is non-connected.

Basically, the same example can be used for point sets by adding a new vertex $p_7$ on $s$. To complete a triangulation, join $p_7$ to all the other vertices to obtain a set $S$. By construction, it is not possible to perform flips in any of the quadrilaterals having $p_7$ as a
vertex (the new diagonals are outside the quadrilaterals). So, $S$ has the two different triangulations of the original polygon $P$ and no flip is possible in any of those triangulations.
\qquad \end{proof}

Using the previous lemma,
the same reasoning can be extended to the remaining closed connected surfaces by using the fact
that every closed connected surface is topologically equivalent to a sphere, or a connected sum
of tori (handles), or a connected sum of projective planes.

\begin{theorem}
\label{noconne} Any closed connected surface $S$ admits a metric that allows polygons and point sets
whose graphs of (metrical) triangulations are non--connected.
\end{theorem}

\begin{figure}[htb]
\begin{center}
\includegraphics[height=4cm]{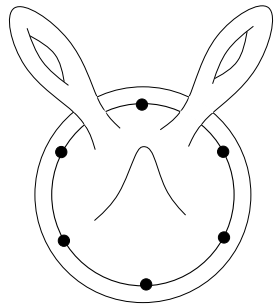}
\caption{The  construction of Fig.~\ref{sph-tor2} on the sphere with two handles.} \label{sph-tor2}
\end{center}
\end{figure}

\begin{proof}
We can modify the surface $M$ described in the proof of Lemma~\ref{lsnc}, by adding to it as many handles or projective planes as needed in order to obtain a surface homeomorphic to $S$.

By virtue of this fact, and mimicking the argument we followed on the sphere, it is possible to find a
metric on each closed and connected surface that allows polygons and point sets with non-connected
graph of triangulations (see Figure~\ref{sph-tor2}).
\qquad \end{proof}

Note that the reasoning used on the proof of Theorem~\ref{noconne} can be easily extended to any kind of surface.

\section{Connectivity of the graph of triangulations on locally Euclidean surfaces}
\label{locally}
As it has been said in the Introduction,
some of the most common and useful surfaces are the locally Euclidean surfaces, because of the
advantage of their planar representations. It is interesting to emphasize the different behavior
that these surfaces show when we study the graph of triangulations of a polygon: while it is
connected both in the cylinder and in the flat torus, polygons with non-connected graph can easily
be constructed in the two non-orientable surfaces. The behavior of the graph of triangulations in the torus is remarkable since, that graph  is connected for polygons with the metric of the flat torus, but this is not longer true for the skew torus. On
the other hand, the graph of a point set is non-connected in general  but, as we shall see next, we can describe all the connected components in the case of the cylinder.

\subsection{The cylinder} \label{cylinder}

Let $\vec{a}$ be the vector that generates the cylinder. Given an orthogonal reference system being the $OX$ axis parallel to $\vec{a}$, a geodesic arc is a segment if and only if its vertical projection is smaller than $|\vec{a}|/2$. In order to add a new diagonal to a triangulation, a procedure to determine if the geodesic arc joining two vertices it is a segment is to check if its vertical projection is contained inside the vertical projection of a previously existing diagonal (and, therefore, a segment).

\subsubsection{Polygons}

If the planar copies of a polygon $P$ on the cylinder are (each of them) strictly contained in vertical bands of length $|\vec{a}|/2$, then any
internal diagonal is admissible and planar arguments can straightforwardly be
used to establish the connectivity of the graph of triangulations~\cite{albertoyclara}.

However, although many different proofs are known for planar polygons in the plane, the authors are not aware of any proof that can be adapted  for the general case. Actually, it is not even obvious that in this
general situation a polygon can always be triangulated, although, in this case, essentially the same ideas as in
the plane provide a proof of this fact.

\begin{lemma} \label{cil-triang}
Any Euclidean polygon of $n\geq 4$ vertices on the cylinder has an admissible diagonal. Hence, any
Euclidean polygon on the cylinder is triangulable.
\end{lemma}

\begin{figure}[htb]
\begin{center}
\includegraphics[height=4.8cm]{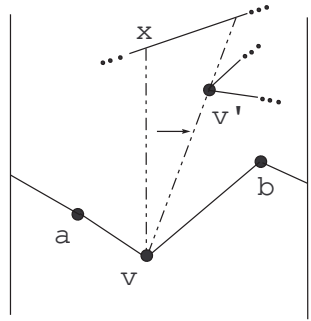}
\caption{$vv'$ is a diagonal of the polygon.} \label{rota}
\end{center}
\end{figure}

\begin{proof}
This proof is based on the proof of Meister's Lemma~\cite{meister}, which establishes the same
result for simple polygons in the plane.

Consider a Euclidean polygon $P$ already developed into the plane.
Let $v$ be a convex vertex such that the two edges incident on it go
upwards (recall that a vertex is {\em convex} if its interior angle
is less than $\pi$ radians; otherwise, the vertex is {\em reflex}) .
Let $a$ and $b$ be the vertices adjacent to $v$.
(Figure~\ref{diag-cyl}).

\begin{figure}[htb]
\begin{center}
\includegraphics[height=4cm]{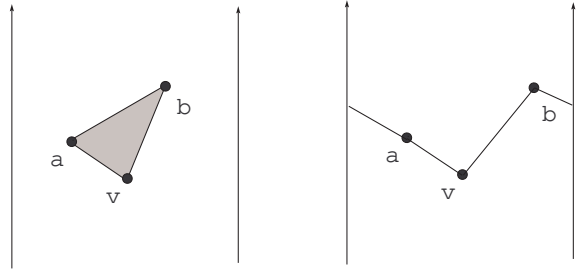}
\caption{The segment matching $a$ and $b$ may determine or not a
bounded triangle.} \label{diag-cyl}
\end{center}
\end{figure}

If $\overline{ab}$ is an admissible diagonal (a segment contained in $P$), we have
finished. Otherwise, either $\overline{ab}$ intersects $\partial P$, or
it is exterior to $P$.

If $\overline{ab}$ intersects $\partial P$, the argument given in~\cite{meister} can be mimicked:
Start sweeping a line from $v$, keeping it parallel to the line through $\overline{ab}$ until it
reaches another vertex $x$ of $P$ (it must exist since $P$ has at least four vertices). Then,
$\overline{vx}$ is an admissible diagonal (Figure~\ref{inn-diag}).

\begin{figure}[htb]
\begin{center}
\includegraphics[height=4cm]{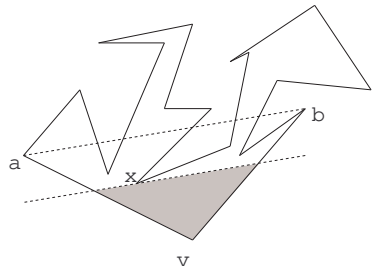}
\caption{$vx$ is an inner diagonal of the polygon.}
\label{inn-diag}
\end{center}
\end{figure}

If $\overline{ab}$ is exterior to $P$, consider the vertical {\em ray} (half-line) with $v$ as
endpoint and let $x$ be the first point of the boundary of $P$ that it reaches. If $x$ is a vertex
then $\overline{vx}$ is an admissible diagonal. Otherwise, rotate the ray either to the right or
to the left until it intersects another vertex $v'$ of $P$ (Figure~\ref{rota}). The vertical projection of $\overline{vv'}$ is contained inside the vertical projection of the diagonal containing $x$, so $\overline{vv'}$ is an
admissible diagonal.
\qquad \end{proof}

The connectivity of the graph of triangulations of a Euclidean
polygon on the cylinder is established by the next theorem. As in
the plane, three consecutive vertices of the polygon $a$, $v$, $b$
are said to form an {\em ear} if $\overline{ab}$ is an admissible
diagonal: $v$ is called the ear {\em tip}. Two ears are {\em
non-overlapping} if their triangle interiors are disjoint.

\begin{theorem}\label{res14tri}
Any two triangulations of a Euclidean polygon of $n\geq 4$ vertices
on the cylinder can be transformed one into each other by a sequence
of flips.
\end{theorem}

\begin{proof}
Observe first that any triangulation of a Euclidean polygon of $n\geq 4$ vertices has at least
two non-overlapping ears. The proof of this fact in the plane \cite{meister} is essentially
topological and works the same in the planar representation of the polygon, no matter the metric.

Consider two triangulations $T_1$ and $T_2$ of a Euclidean polygon $P$ with $n\geq 4$ vertices.
We proceed by induction on the number of vertices (the base case is obvious). The inductive hypothesis
straightforwardly leads to the result if $T_1$ and $T_2$ have a diagonal in common, since this
diagonal divides $P$ into two smaller polygons.

Suppose then that $T_1$ and $T_2$ do not share any diagonal. Let $v_1$ (resp., $v_2$) be the tip
of an ear in $T_1$ (resp., in $T_2$). Since both $T_1$ and $T_2$ have at least two non-overlapping
ears, $v_1$ and $v_2$ can be assumed to be non-adjacent.

Consider $T_2-\{v_2\}$ that is the triangulation of a polygon of $n-1$ vertices. By induction,
$T_2-\{v_2\}$ can be transformed by flips into another triangulation $T'_2-\{v_2\}$ having an ear
in $v_1$ (Lemma~\ref{cil-triang} assures that such a triangulation exists). As a consequence, $T_2$ can be transformed into another triangulation $T'_2$ having an
ear in $v_1$. Now $v_1$ is the tip of an ear both in $T_1$ and in $T'_2$ what means that they
share a diagonal.

It follows from the induction hypothesis that it is possible to go from $T_1$ to $T_2$ by flips through $T'_2$:
$T_1\hookrightarrow T'_2 \hookrightarrow T_2$.
\qquad \end{proof}

It is easily seen that the proof of Theorem \ref{res14tri} implies an $O(2^n)$ upper bound for the
number of flips. This bound is far to be tight on the cylinder since we have the same example as in the
plane~\cite{hurt2} for an $\Omega(n^2)$ bound  and any simple planar polygon, conveniently reduced, can be embedded in half a
cylinder with the same (global) metric~\cite{albertoyclara}.

\subsubsection{Point sets}

Regarding the connectivity of the graph of triangulation of a point set on the cylinder, three
different situations can be presented.

Consider the circle determined by cutting the cylinder with a plane orthogonal to its axis. If the
smallest arc covering the orthogonal projection of the points on that circle is smaller than
$\pi$, then the set is in {\em Euclidean position} and it has a planar behavior
\cite{euclideanposition}. So, their graph of triangulation is connected.

\begin{figure}[htb]
\begin{center}
\includegraphics[height=4.8cm]{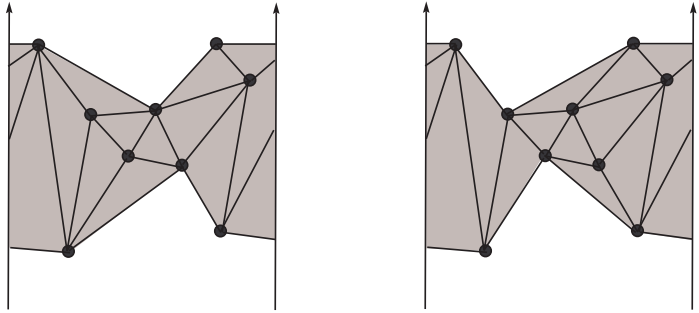}
\caption{The points are the same, but the triangulated region is different. As there can not be
performed any flip on the segments of the boundary, the graph of triangulations is disconnected.}
\label{tcilin}
\end{center}
\end{figure}

If the set is not in Euclidean position, then there exist three of its points such that the
polygonal line joining then wraps around the cylinder. We call this an {\em essential polygonal line}. The
triangulated region of a set that is not in Euclidean position is bounded by two closed essential
polygonal lines. The polygonal lines bounding the triangulated region are not uniquely determined by the
points, and different triangulations of the same set may be bounded by different polygonal lines,
describing different regions, as shown in Figure~\ref{tcilin}. In this case, as it is not
possible to perform flips over the segments of the boundary, the graph of triangulations of the
set is non--connected.

Finally, the third case occurs when we consider triangulations of a point set that share the same
boundaries. Then it is possible to carry one into each other by a sequence of flips as it is established
by the following theorem.

\begin{theorem} \label{cil-nube}
Given two triangulations of a point set on the cylinder, with the same boundaries, it is possible
to transform one into the other by a sequence of flips.
\end{theorem}

\begin{figure}[htb]
\begin{center}
\includegraphics[height=2.7cm]{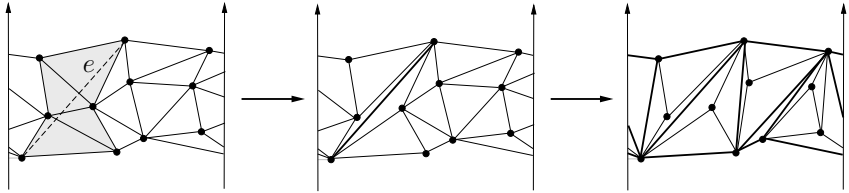}
\caption{Construction of $T_D^1$.} \label{e}
\end{center}
\end{figure}

\begin{proof}
Let $T_1$ y $T_2$ be two triangulations of a point set  $S$, and let
$U$ and $L$  be the (upper and lower respectively) essential
polygonal lines bounding $S$. Denote $D$  the region bounded by $U$ and $L$.

The first step in the proof is to triangulate $D$ using only admissible diagonals from $U$ to $L$, or between two vertices of $U$ or two vertices of $L$, without considering interior points of $S$. In order to get this goal, we use the same arguments as in  Lemma~\ref{cil-triang}. Denote $T_D$ the triangulation so obtained.

Now, let $e$ be a diagonal in  $T_D$, and let $P_e$ be the polygon defined by the union of all the triangles in $T_1$ that intersect $e$. We can extend $e$ to a triangulation of $P_e$ by using Lemma~\ref{cil-triang}. Obviously, by Theorem~\ref{res14tri}, we can transform one of those triangulations into the other, and so we can transform $T_1$ into another triangulation containing  the diagonal $e$. We can do the same process for all diagonals in $T_D$  to obtain a triangulation $T_D^1$, and the same starting from $T_2$ to obtain a new triangulation $T_D^2$ such that it is possible to transform  $T_i$ into $T_D^i$ using admissible flips ($i=1,2$)
(Figure~\ref{e}). The only differences between $T_D^1$ and $T_D^2$ are in diagonals that are contained into triangles of $T_D$. But the vertical projection of any diagonal inside a triangle of $T_D$ is contained on the vertical projection of one of the sides of the triangle, so all the diagonals are admissible, the we are in the same situation that in the plane and, therefore, we can transform $T_D^1$ into $T_D^2$ using flips.
\qquad \end{proof}

Of course, Theorem~\ref{cil-nube} provides a method to characterize the connected components of
the graph of triangulations of a point set on the cylinder in terms of the number of upper and
lower polygonal lines; in \cite{polardiagram} it is given a tool, the polar diagram, that allows to
count those chains.

\subsection{The Torus} \label{torus}

Throughout this section, unless otherwise stated, we assume a flat
torus generated by a pair $\vec{a}$ and $\vec{b}$ of orthogonal
vectors. Those vectors are considered to be horizontal and vertical,
respectively. A fundamental domain for this surface is an
isothetic rectangle of dimension $|\vec{a}|\times|\vec{b}|$. Recall that an arc of geodesic is a segment if and only if its vertical projection is smaller than $|\vec{a}|/2$ and its horizontal projection is smaller than
$|\vec{b}|/2$. Note that this implies that if a planar copy of a set on the torus is contained on an isothetic rectangle of dimension
$|\vec{a}|/2\times|\vec{b}|/2$ then any diagonal is admissible and the set has a planar behavior. Analogously, if a planar copy of the set is inside a vertical (resp., horizontal) strip of width $|\vec{a}|/2$ (resp., $|\vec{b}|/2$) the behavior is equivalent to be on the cylinder~\cite{euclideanposition,albertoyclara}.

We will center our efforts in proving the connectivity of the graph of triangulations of polygons on
the flat torus. It is still an open problem whether the graph of triangulations is connected or not when we consider point sets instead of polygons.
It is worthy to note that if the torus is generated by non-orthogonal vectors (skew torus), sets and
polygons with disconnected graph of triangulations appear.

\begin{figure}[htb]
\begin{center}
\includegraphics[height=6cm]{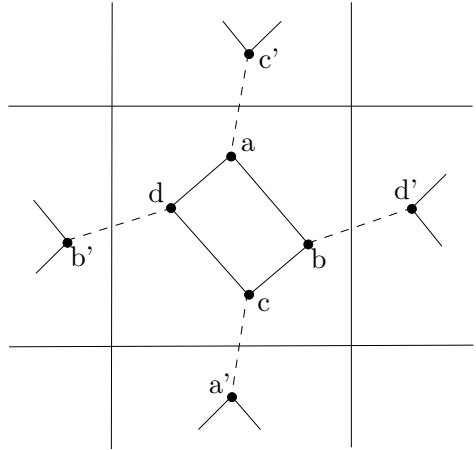}
\caption{There is not any admissible diagonal in this
quadrilateral.} \label{ntpolyg}
\end{center}
\end{figure}

To establish the connectivity of the graph of triangulations of a polygon on the flat torus is
more complicated than on the cylinder. The most common proofs of the connectivity of the graph of
triangulations (either in the plane or on the cylinder) are based, in some sense, in the fact that
every polygon is triangulable, but this is no longer true on the torus (as
Figure~\ref{ntpolyg} shows), and inductive reasonings fail. To be more precise, a maximal set of
admissible diagonals does not necessarily divide the interior of the polygon into triangular
regions. That is why the term {\em triangulation of a polygon\/} was defined in
\S\ref{trianpolygon} as a partition of the interior of the polygon into triangular regions
rather than maximal sets of admissible diagonals.

\begin{figure}[htb]
\begin{center}
\includegraphics[height=4.8cm]{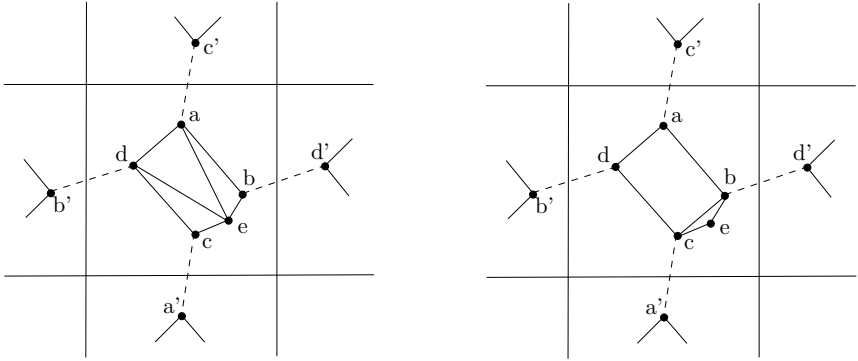}
\caption{This pentagon is triangulable, but there exists no
triangulation containing the diagonal $\overline{cd}$.}
\label{determ}
\end{center}
\end{figure}

\begin{figure}[htbp]
\begin{center}
\includegraphics[height=4.8cm]{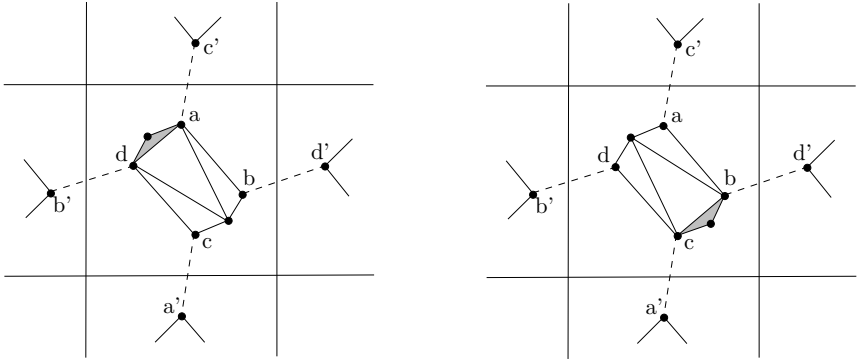}
\caption{If we fix the selected ear in the first triangulation, as
we did on the cylinder, it is not possible to perform any flip in
the rest of the triangulation to carry it to one that contains the
selected ear in the second triangulation.} \label{contraej}
\end{center}
\end{figure}

Moreover, even if the polygon admits a triangulation, some
admissible diagonals may not take part of any of them, as it is
shown in Figure~\ref{determ} where the polygon is
triangulable but the admissible diagonal $\overline{cd}$ does not
participate in any triangulation. In Figure~\ref{contraej} we can
see two different triangulations of a polygon where the method
used on the cylinder fails. Since there is no other admissible
diagonal than those shown in the figure, it is not possible to
transform one of the  triangulations into the other by keeping one ear
in common with one of the two triangulations along the process.

This anomalous behavior of the torus restricts ourselves to the study of the graph of
triangulations of triangulable polygons and forces the search for new techniques to establish the
connectivity of the graph of triangulations.

However, it is possible to give necessary conditions for a polygon to be triangulable on the flat
torus. We define a {\em quadrant} in the torus as an isothetic rectangle of dimension
$|\vec{a}|/2\times|\vec{b}|/2$, with $|\vec{a}|$ and $|\vec{b}|$ being the generating vectors of
the flat torus. It is easily seen that a diagonal of a polygon on the flat torus is admissible if
and only if it fits in a quadrant. Then the following result is straightforward.

\begin{proposition}
If a polygon $P$ on the torus contains the center of an empty (with no vertex of the polygon inside it) quadrant then $P$
is not triangulable.
\end{proposition}

\begin{proof}
Suppose, on the contrary to our claim, that $P$ is a triangulable polygon containing the center
$O$ of an empty quadrant. Then there must exist a triangle $T_O$ of the triangulation of $P$ in
which $O$ lies in. Now it can easily be checked  that some of the edges of $T_O$ are not
admissible, conclusion contrary to our assumption.
\qquad \end{proof}

We now introduce some definitions and preliminary results that
will lead to a proof of the connectivity of the graph of triangulations of
a triangulable polygon on the flat torus.

Let $P$ be a Euclidean polygon on the flat torus. We consider a copy of $P$ in the plane with the usual reference system.  A convex vertex $v$ of $P$ is said to be a {\em top} vertex (resp., {\em bottom}, {\em right} and {\em left})
if the two edges incident on it go downwards (resp., upwards, leftwards and rightwards).
Top, bottom, right and left vertices will be called {\em extreme} vertices.

A vertex $u$ of $P$ is said to be {\em earable\/} if the segment
joining its  two adjacent vertices is an admissible diagonal.

\begin{figure}[htb]
\begin{center}
\includegraphics[height=3.5cm]{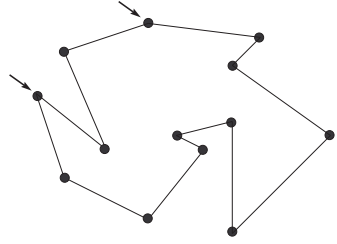}
\caption{A polygon with the {\em top\/} vertices marked.} \label{top}
\end{center}
\end{figure}

\begin{lemma}
\label{l1tp} Any triangulable quadrilateral on the flat torus admits a
triangulation having an ear in one of its extreme vertices.
\end{lemma}

\begin{figure}[htb]
\begin{center}
\includegraphics[height=2.5cm]{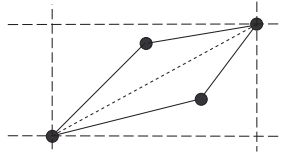}
\caption{Any triangulable quadrilateral has an extreme earable vertex.} \label{cuad}
\end{center}
\end{figure}

\begin{proof}
Since the two ears in a triangulation of a quadrilateral are at opposite vertices, and since every
polygon has at least two extreme vertices, the only non-trivial case is that of a quadrilateral
with two non-extreme vertices opposite to one another. Such a quadrilateral is contained in the
isothetic rectangle having the two extreme vertices as corners (Figure~\ref{cuad}).
This implies that if the diagonal joining the extreme vertices is admissible, then the diagonal joining the non-extreme vertices is admissible as well. One of them must be admissible because the quadrilateral is assumed to be triangulable.
\qquad \end{proof}

Now, we extend the result given in Lemma~\ref{l1tp} to a general
polygon.

\begin{lemma}
\label{l2tp} Any triangulable polygon on the flat torus admits a
triangulation having an ear in one of its extreme vertices.
\end{lemma}

\begin{proof}
We proceed by induction in the number of vertices. The base case is given by the previous lemma,
so consider a triangulable polygon $P_{n+1}$ with $n>4$ vertices and let $T$ be a triangulation of
$P_{n+1}$. Let $v$ be the tip of an ear of $T$ and let $a$ and $b$ be its adjacent vertices. If $v$ is an
extreme vertex in $P_{n+1}$, we finish. Otherwise, remove $v$ and its incident segments $\overline{av}$ and $\overline{vb}$ from $P_{n+1}$ to obtain a polygon
$P_n$ with $n$ vertices. By the induction hypothesis, $P_n$ admits a triangulation $T'$ having an
ear in one extreme vertex $v'$.

If $v'$ is other than $a$ and $b$, then $v'$ is also an extreme
vertex in $P_{n+1}$ and $T''=T'\cup \{\triangle avb\}$ is a
triangulation of $P_{n+1}$ having an ear in $v'$, so the result
holds.
Therefore, assume $T'$ has an extreme ear $\triangle cab$ at $a$.

Without loss of generality, suppose $a$ is a right vertex of
$P_n$.  Let $R$ be the lower triangle defined by the edge $\overline{ab}$ in
the isothetic rectangle having $\overline{ab}$ as one of its diagonals ($R$
is the shaded area in Figure~\ref{nuevo-v}). Since $v$ is not
extreme, it must be inside $R$. But this implies that $a$ is also
a right vertex in $P_{n+1}$, as it is shown in Figure~\ref{nuevo-v}.

\begin{figure}[htb]
\begin{center}
\includegraphics[height=3.5cm]{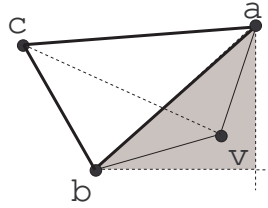}
\caption{If $v$ is outside $R$ it is an extreme vertex. Otherwise, a
flip can be performed to obtain an extreme ear in $a$.}
\label{nuevo-v}
\end{center}
\end{figure}

The vertical (resp., horizontal) projection of $\overline{vc}$ is
inside the vertical (resp., horizontal) projection of $\overline{ac}$ (resp., $\overline{ab}$), and since
both $\overline{ac}$ and $\overline{ab}$ are admissible diagonals, so $\overline{vc}$ is admissible, and $\overline{ab}$ can
be flipped in $T''=T'\cup \{\triangle avb\}$ in order to get the
extreme ear $\triangle cav$ in the triangulation $T''$ of $P_{n+1}$.
\qquad \end{proof}

Now, we can prove the key result in order to obtain the connectivity of the graph of
triangulations of a polygon on the flat torus.

\begin{lemma}
\label{l3tp} Let $P$ be a triangulable polygon on the flat torus
having an extreme earable vertex $u$. Then, any triangulation of $P$
can be transformed by a sequence of flips into a triangulation
having an ear in $u$.
\end{lemma}

\begin{figure}[htb]
\begin{center}
\includegraphics[height=3.5cm]{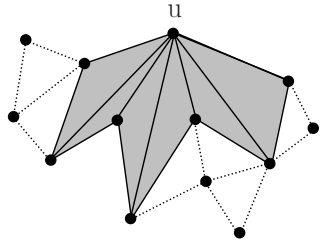}
\caption{$P'$ is made by the triangles incident to $u$.} \label{u3}
\end{center}
\end{figure}

\begin{proof}
For the sake of simplicity we assume that $u$ is a top vertex. Let
$T$ be a triangulation of $P$ such that $u$ is not an ear in $T$. Consider the subpolygon $P'$ of
$P$ covered by the triangles of $T$ incident to $u$ (Figure~\ref{u3}). It is clear that $u$ is the topmost vertex
of this subpolygon, and since all the vertices of $P'$ are joined to $u$ by admissible diagonals, $P'$ is contained in a horizontal strip of width $|\vec{b}|/2$. Thus, we can consider $P'$ as it is embedded on a cylinder generated by $\vec{a}$ and, by our results in \S\ref{cylinder}, $P'$ admits a triangulation which has an ear at the earable vertex $u$ and this triangulation is connected by
flips to the restriction of $T$ to $P'$. This finishes the proof.
\qquad \end{proof}

Now, we are in the condition to enunciate the main result of this section.

\begin{theorem}
\label{t1tp} The graph of triangulations of a polygon on the flat torus
is either empty or connected.
\end{theorem}

\begin{proof}
Let $T_1$ and $T_2$ be two triangulations of a polygon $P$ on the flat torus. Let $u$ be an
extreme earable vertex in $P$, which exists by Lemma~\ref{l2tp}. By virtue of Lemma~\ref{l3tp},
$T_1$ (resp., $T_2$) can be transformed by a sequence of flips into another triangulation $T'_1$
(resp., $T'_2$) having an ear in $u$. Therefore $T'_1$ and $T'_2$ are connected using flips by the
inductive hypothesis.
\qquad \end{proof}

Regarding the connectivity of the graph of triangulations of a point set $S$ in the flat torus there
are three possible situations:
    \begin{remunerate}
    \item If $S$ is inside a quadrant. Then $S$ is in Euclidean position and it has a planar
    behavior \cite{euclideanposition}, so the graph is connected.
    \item If a planar copy of $S$ is inside a vertical (resp., horizontal) strip of width $|\vec{a}|/2$ (resp., $|\vec{b}|/2$) the situation is equivalent to the cylinder. The graph is connected if and only if the borders of the triangulated region are fixed (\S\ref{cylinder}).
    \item In other case the connectivity of the graph of triangulations is still an open problem. Our conjecture is that this graph is connected.
    \end{remunerate}

\begin{figure}[htb]
\begin{center}
\includegraphics[height=7.6cm]{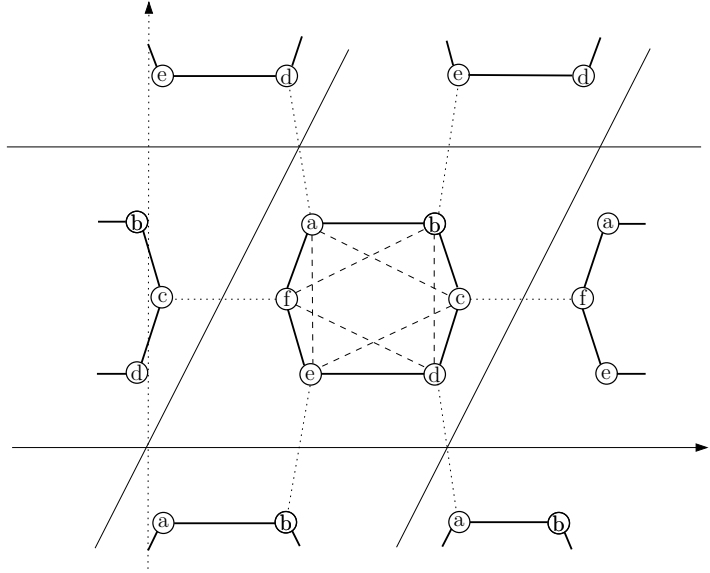} \caption{The previous hexagon with
all the possible segments between its vertices.} \label{sk-tor}
\end{center}
\end{figure}

Nevertheless, as we pointed out in \S\ref{ccsurfaces}, the connectivity of the graph of
triangulations is not preserved if the torus is generated by two non-orthogonal translations.
In this way, consider the planar representation of a skew torus generated by two translations with
vectors forming an angle of $\arccos\sqrt{\frac{1}{5}}$. In order to simplify the coordinates of
the vertices, we choose an horizontal unitary vector and the other one with modulo
$\frac{2\sqrt{5}}{5}$ and hence the height of a fundamental region is one unit. Using the usual
reference system we can draw an hexagon of vertices $a(\frac{1}{2}+\varepsilon,\frac{3}{4})$
$b(1-\varepsilon,\frac{3}{4})$, $c(1+\frac{\varepsilon}{3},\frac{1}{2})$,
$d(1-\varepsilon,\frac{1}{4})$, $e(\frac{1}{2}+\varepsilon,\frac{1}{4})$, and
$f(\frac{1}{2}-\frac{\varepsilon}{3},\frac{1}{2})$, with $\varepsilon < \frac{1}{8}$. Since the
diagonals $ad$, $be$ and $cf$ are not admissible, it is not possible to perform flips in any
of the two triangulations depicted in Figure~\ref{sk-tor}.

\begin{figure}[htb]
\begin{center}
\includegraphics[height=7.6cm]{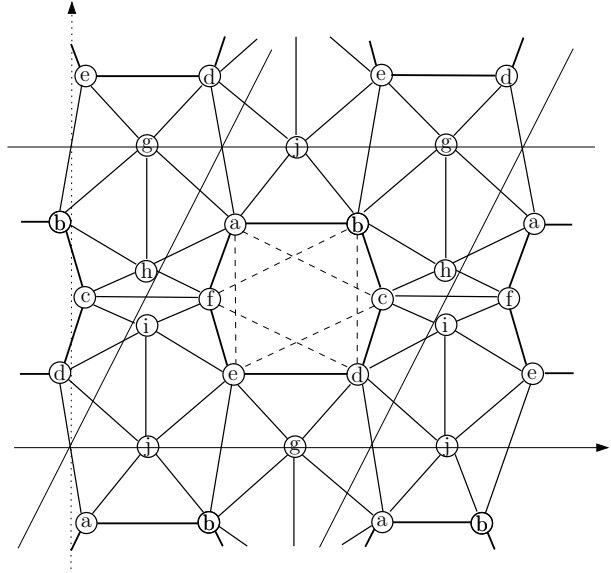} \caption{It is not possible to carry one of
the triangulations of the central polygon into the other by flips.} \label{tnb-tort}
\end{center}
\end{figure}

And, likewise we have done in \S\ref{ccsurfaces}, new points can be added to the previous
construction to obtain a point set with a non-connected graph of triangulations. We include points
$g(0,\frac{3}{4})$, $h(\frac{1}{4},\frac{3}{5})$, $i(\frac{1}{4},\frac{2}{5})$ and
$j(0,\frac{1}{4})$ as it is shown in Figure~\ref{tnb-tort}. The central hexagon (bold lines)
still admits only six diagonals, giving rise to only two different triangulations, and the segments
of the boundary of the hexagon cannot be flipped. So the graph of triangulations of the set has two connected components.

Therefore, the previous example shows (applying a suitable angle
transformation if necessary) the following result.

\begin{theorem} \label{skewtorus}
It is possible to find a polygon and a point set on a skew torus such that their graphs of
(metrical) triangulations are non-connected.
\end{theorem}

It is worthy to point out that the previous result leads  to another proof of
Theorem~\ref{noconne}.

\subsection{Non-orientable locally Euclidean surfaces}
\label{non-orient}

It is easy to embed in the twisted cylinder and in the Klein bottle a polygon which graph of triangulations is non-connected. We look for a situation similar to the one used previously for the skew torus (Figure~\ref{sk-tor}), an hexagon of vertices (in clockwise order) $a,b,\dots,f$ such that all the diagonals are admissible but the diagonals $(a,d)$, $(b,e)$ and $(c,f)$. This hexagon has only two triangulations and it is not possible to perform any flip.

\begin{figure}[htb]
\begin{center}
\includegraphics[height=4cm]{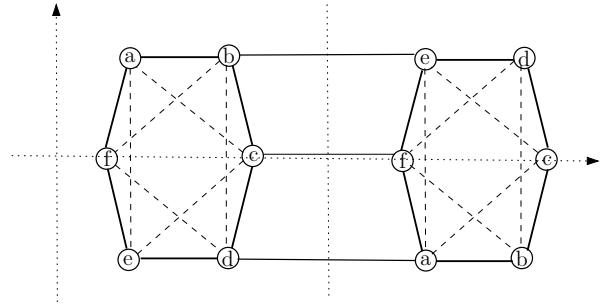}
\caption{The segments $ad$, $be$ and $cf$ are outside the polygon.}
\label{pol-kl}
\end{center}
\end{figure}

Consider the twisted cylinder
with the usual coordinate system, and assume a glide reflection with unitary vector. The hexagon
of vertices $a(\frac{1}{4}+\varepsilon,\frac{1}{4})$, $b(\frac{3}{4}-\varepsilon,\frac{1}{4})$,
$c(\frac{3}{4}+\frac{\varepsilon}{3},0)$, $d(\frac{3}{4}-\varepsilon,-\frac{1}{4})$,
$e(\frac{3}{4}+\varepsilon,-\frac{1}{4})$, and $f(\frac{1}{4}-\frac{\varepsilon}{3},0)$, with
$\varepsilon < \frac{1}{16}$ has a non-connected graph of triangulations (Figure~\ref{pol-kl}).

The same construction can be easily obtained on the Klein bottle. A similar study can be
extended to the other surfaces obtained as the quotient of the plane over a group of
motions ({\em Euclidean 2-orbifolds\/}), if the group contains a glide reflection. In particular, this hexagon can also be embedded on
the projective plane with the quotient metric.

\begin{figure}[htb]
\begin{center}
\includegraphics[height=4cm]{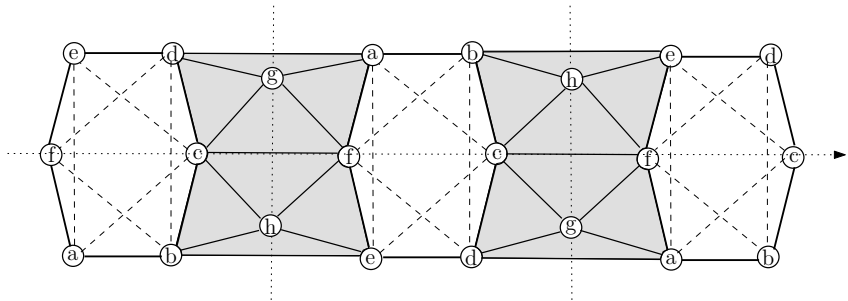} \caption{The only flips allowed are restricted to the
shaded regions.} \label{tnb-cilr}
\end{center}
\end{figure}

By adding only two new points, as Figure~\ref{tnb-cilr} shows, the previous example can
be extended to a point set with a non-connected graph of triangulations. Again, as we saw in
\S\ref{torus} for the skew torus, the central hexagon has two possible triangulations and
no flip is possible inside it. And, since the segments of the boundary of the hexagon cannot be
flipped, the two triangulations belong to different connected components.

\begin{figure}[htb]
\begin{center}
\includegraphics[height=8cm]{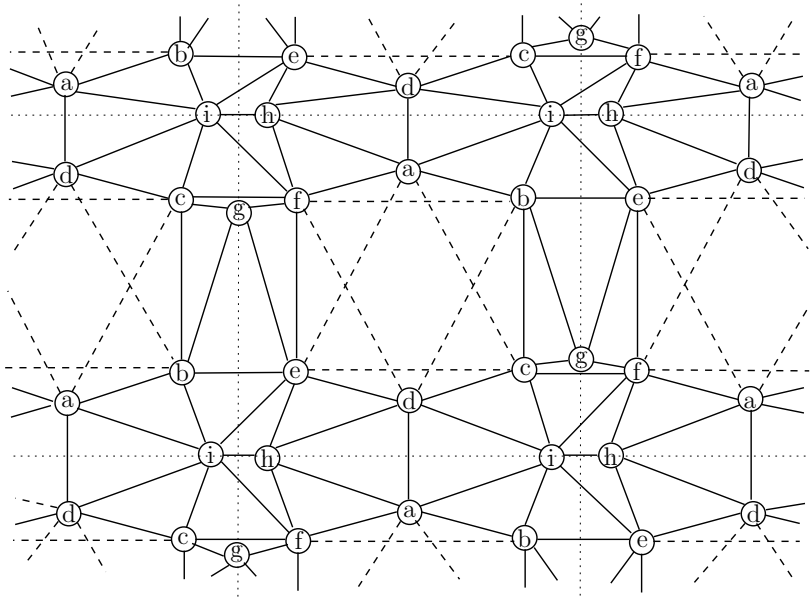} \caption{It is not possible to carry one of
the triangulations of the central polygon into the other by flips.} \label{tnb-kl}
\end{center}
\end{figure}

Slightly more complicated is the example given for the Klein bottle (Figure~\ref{tnb-kl}),
but the reasoning is the same; the segments that form the central hexagon cannot be flipped, hence the set has two disjoint triangulations.

\section{Conclusions and open problems}
\label{open-p}

In this paper we have studied the connectivity of the graph of triangulations of polygons and
point sets on a surface. We have seen that, in general, this graph is non-connected. More
precisely, we have proven that any surface admits a metric such that there exits a polygon and a
point set on that surface (and with that metric) with non-connected graph of triangulations. There
exist some remarkable exceptions. First of all, of course, the plane and the sphere, and then
polygons on the cylinder and the flat torus, these surfaces being the only ones in which Lawson's
method~\cite{lawson2} to obtain an optimal triangulation could be applied. Nevertheless it is not clear the practical applications of this method since we have not found reasonable bounds for the diameter of the graph of triangulations in these surfaces.

Some problems are still unsolved. The main question that remains after Theorem~\ref{noconne} is whether
is possible or not to define a metric in a surface forcing the graph of triangulations of any
polygon to be connected, as it has been shown for the torus.

And, back to the flat torus, the connectivity of the graph is not established if it is associated
to triangulations of point sets instead of polygons.

\bibliographystyle{siam}

\begin{thebibliography}{10}

\bibitem{avis-fukuda}
{\sc D.~Avis and K.~Fukuda}, {\em Reverse search for enumeration}, Discrete
  Applied Math, 6 (1996), pp.~21--46.

\bibitem{bern}
{\sc M.~Bern and D.~Eppstein}, {\em Mesh generation and optimal triangulation},
  in Computing in Euclidean Geometry, D. Z. Du and F. K. Hwang, World
  Scientific, 1992, pp.~23--90.

\bibitem{euclideanposition}
{\sc C.~Cortes, A.~Marquez, and J.~Valenzuela}, {\em Euclidean position in
  euclidean 2-orbifolds}, Computational Geometry, 27 (2004), pp.~27--41.

\bibitem{docarmo2}
{\sc M.~P. do~Carmo}, {\em Differential geometry of curves and surfaces},
  Prentice, inc. USA, 1976.

\bibitem{hurt2}
{\sc M.~Noy F.~Hurtado and J.~Urrutia}, {\em Flipping edges in triangulations},
  Discrete and Computational Geometry, 22 (1999), pp.~333--346.

\bibitem{fortune2}
{\sc S.~Fortune}, {\em Voronoi diagrams and {D}elaunay triangulations}, in
  Computing in Euclidean Geometry, D. Z. Du and F. K. Hwang, World Scientific,
  1992, pp.~193--234.

\bibitem{gelatt}
{\sc C.~D. Gelatt, S.~Kirkpatrick, and M.~P. Vecchi}, {\em Optimization by
  simulated annealing}, Science, 220 (1983), pp.~671--680.

\bibitem{albertoyclara}
{\sc C.~I. Grima and A.~M\'{a}rquez}, {\em Computational Geometry on Surfaces},
  Kluwer Academic Publisher, 2001.

\bibitem{polardiagram}
{\sc C.~I. Grima, A.~M\'{a}rquez, and L.~Ortega}, {\em A new 2d tessellation
  for angle problems: the polar diagram}, Computational Geometry, 34 (2006),
  pp.~58--74.

\bibitem{hurt}
{\sc F.~Hurtado and M.~Noy}, {\em Graph of triangulations of a convex polygon
  and tree of triangulations}, Computational Geometry, 13 (1999), pp.~179--188.

\bibitem{vanlaar}
{\sc P.~J. M.~Van Laarhoven and E.~H.~L. Aarts}, {\em Simulated Annealing:
  Theory and Practice}, Kluwer Academic Publ., 1987.

\bibitem{lawson2}
{\sc C.~L. Lawson}, {\em Transforming triangulations}, Discrete Mathematics, 3
  (1972), pp.~365--372.

\bibitem{mazon}
{\sc M.~Maz\'{o}n and T.~Recio}, {\em Voronoi diagrams on orbifolds},
  Computational Geometry, 8 (1997), pp.~219--230.

\bibitem{meister}
{\sc G.~H. Meister}, {\em Polygons have ears}, Amer. Math. Monthly, 82 (1975),
  pp.~648--651.

\bibitem{nikulin}
{\sc V.~V. Nikulin and I.~R. Shafarevich}, {\em Geometries and Groups},
  Springer Series in Soviet Mathematics, Springer, Berlin, 1987.

\bibitem{okabe}
{\sc A.~Okabe, B.~Boots, and K.~Sugihara}, {\em Spatial Tesselations. Concepts
  and Applications of Voronoi Diagrams}, John Wiley \& Sons, 1992.

\bibitem{pocch}
{\sc M.~Pocchiola and G.~Vegter}, {\em Computing the visibility graph via
  pseudo-triangulations}, in Proceedings of the 11th. ACM Symp. on Comp.
  Geometry, 1996, pp.~248--257.

\bibitem{santos:bistellar}
{\sc F.~Santos}, {\em Geometric bistellar flips. the setting, the context and a
  construction}, in International Congress of Mathematicians, Marta
  Sanz-Sol\'e, Javier Soria, Juan~Luis Varona, and Joan Verdera, eds.,
  vol.~III, European Mathematical Society, 2006, pp.~931--962.

\bibitem{sleator}
{\sc D.~D. Sleator, R.~E. Tarjan, and W.~P. Thurston}, {\em Rotations distance,
  triangulations and hyperbolic geometry}, J. Am. Math. Soc., 1 (1988),
  pp.~647--682.

\end{thebibliography}

\end{document}